\newcommand{\ga}{\gamma}                     \newcommand{\Ga}{\Gamma}
\newcommand{\de}{\delta}                     
\newcommand{\lb}{\lambda}

              \newcommand{\vphi}{\varphi}

\newcommand{\cal}{\mathcal}

                 \newcommand{\calv}{{\cal V}}

             \newcommand{\Dom}{{\rm Dom}}

               \newcommand{\incl}{\subseteq}
\newcommand{\es}{\emptyset}                  \newcommand{\sm}{\setminus}
              \newcommand{\limpl}{\Longrightarrow}
          \newcommand{\lequi}{\Longleftrightarrow}
\newcommand{\oo}{\infty}

                   \newcommand{\sk}{\smallskip}
             \newcommand{\n}{\noindent}

\newcommand{\barr}{\begin{array}}            \newcommand{\earr}{\end{array}}
\newcommand{\beq}{\begin{equation}}          \newcommand{\eeq}{\end{equation}}
\newcommand{\bit}{\begin{itemize}}           \newcommand{\eit}{\end{itemize}}
\newcommand{\blemma}{\begin{lemma}}          \newcommand{\elemma}{\end{lemma}}
\newcommand{\bproof}{\begin{proof}}          \newcommand{\eproof}{\end{proof}}
\newcommand{\bprop}{\begin{proposition}}     \newcommand{\eprop}{\end{proposition}}
\newcommand{\btab}{\begin{tabular}}          \newcommand{\etab}{\end{tabular}}
\newcommand{\btheorem}{\begin{theorem}}      \newcommand{\etheorem}{\end{theorem}}

\def\Roo{R \cup \{\oo\}}             \def\R+oo{R_+ \cup \{\oo\}}

\def\dtends  {\stackrel {\it d}{\longrightarrow}}
\def\etends  {\stackrel {\it e}{\longrightarrow}}

\def\(V)tends  {\stackrel {(\calv)}{\longrightarrow}}

\documentclass[reqno]{amsart}

\newtheorem{theorem}{\bf Theorem}

\newtheorem{lemma}{\bf Lemma}
\newtheorem{proposition}{\bf Proposition}

\begin{document}

\title
[ALMOST METRIC VERSIONS OF ZVP]
{ALMOST METRIC VERSIONS OF \\
ZHONG'S VARIATIONAL PRINCIPLE}

\author{Mihai Turinici}
\address{
"A. Myller" Mathematical Seminar;
"A. I. Cuza" University;
11, Copou Boulevard; 
700506 Ia\c{s}i, Romania
}
\email{mturi@uaic.ro}




\begin{abstract}
A refinement of  Zhong's variational principle
[Nonlin. Anal., 29 (1997), 1421-1431] is given,
in the realm of almost metric structures.
Applications to equilibrium points are also provided.
\end{abstract}

\maketitle

\section{Introduction}
\setcounter{equation}{0}

Let $M$ be some nonempty set; and $d:M\times M\to R_+$,
a {\it metric} over it (in the usual sense).
Further, take some function  $\vphi:M \to \Roo$ with
\bit
\item[(a01)]
$\vphi$ is inf-proper
($\Dom(\vphi)\ne \es$\ and\  $\vphi_*:=\inf [\vphi(M)]> -\oo$).
\eit
The following 1979 statement in
Ekeland \cite{ekeland-1979}
(referred to as Ekeland's variational principle; in short: EVP)
is our starting point. Assume that
\bit
\item[(a02)]
$d$ is complete (each $d$-Cauchy sequence is $d$-convergent)
\item[(a03)]
$\vphi$ is $d$-lsc
($\liminf_n \vphi(x_n)\ge \vphi(x)$, whenever $x_n\dtends x$).
\eit

\btheorem \label{t1}
Let these conditions hold. Then,

{\bf i)}
for each  $u\in \Dom(\vphi)$ there exists
$v=v(u)\in \Dom(\vphi)$ with
\beq \label{101}
\mbox{
$d(u,v)\le \vphi(u)-\vphi(v)$ (hence $\vphi(u)\ge \vphi(v)$)
}
\eeq
\beq \label{102}
\mbox{
$d(v,x)> \vphi(v)-\vphi(x)$, \quad for all\ $x\in M\sm \{v\}$
}
\eeq

{\bf ii)}
if $u\in \Dom(\vphi)$, $\rho> 0$ fulfill
$\vphi(u)-\vphi_*\le \rho$, then (\ref{101}) gives
\beq \label{103}
\mbox{
($\vphi(u)\ge \vphi(v)$\ and) \quad  $d(u,v)\le \rho$.
}
\eeq
\etheorem

This principle found some basic applications to control and optimization,
generalized differential calculus, critical point theory and global
analysis; we refer to the quoted paper for a survey of these.
So, it cannot be surprising that, soon after its formulation,
many extensions of EVP were proposed.
For example, the (abstract) order one starts from the fact
that, with respect to the (quasi-) order
\bit
\item[(a04)]
($x,y\in M$) $x\le y$ iff $d(x,y)+\vphi(y)\le \vphi(x)$
\eit
the point $v\in M$ appearing in (\ref{102}) is {\it maximal};
so that, Theorem \ref{t1} is nothing but a variant of the
Zorn maximality principle.
The dimensional way of extension refers to the
ambient space ($R$) of $\vphi(M)$ being substituted by a
(topological or not) vector space;
an account of the results in this area is to be found in
the 2003 monograph by
Goepfert, Riahi, Tammer and Z\u{a}linescu
\cite [Ch 3]{goepfert-riahi-tammer-zalinescu-2003}.
Further, the (pseudo) metrical one consists in the conditions
imposed to the ambient metric over $M$ being relaxed.
The basic result in this direction was obtained in 1992 by
Tataru \cite{tataru-1992},
via Ekeland type techniques;
subsequent extensions of it may be found in the 1996 paper by
Kada, Suzuki and Takahashi \cite{kada-suzuki-takahashi-1996}.
Finally, we must add to this list the  "functional" extension of EVP
obtained in 1997 by
Zhong \cite{zhong-1997}
(and referred to as Zhong's variational principle; in short: ZVP).
Take a  function $t\mapsto b(t)$ from $R_+:=[0,\oo[$ to itself, with
the {\it normality} properties
\bit
\item[(a05)]
$b$ is decreasing and $b(R_+)\incl R_+^0:=]0,\oo[$
\item[(a06)]
$B(\oo)=\oo$, \ \ where $B(t)= \int_0^t b(\tau)d\tau, \ t\ge 0$.
\eit

\btheorem \label{t2}
Let $a\in M$ be given, and $u\in \Dom(\vphi)$, $\rho> 0$ 
be taken so as
$\vphi(u)-\vphi_*\le B(d(a,u)+\rho)-B(d(a,u))$.
There exists then $v=v(u)$ in $\Dom(\vphi)$ with
\beq \label{104}
d(a,v)\le d(a,u)+\rho, \quad \vphi(u)\ge \vphi(v)
\eeq
\beq \label{105}
\mbox{
$b(d(a,v))d(v,x)> \vphi(v)-\vphi(x)$, \quad  for each\ \ $x\in M\sm \{v\}$.
}
\eeq
\etheorem

Clearly, ZVP includes (for $b=1$ and $a=u$)
the local version of EVP based upon (\ref{103}).
The relative form of the same, based upon (\ref{101}) 
also holds (but indirectly);
cf. Bao and Khanh \cite{bao-khanh-2003}.
Summing up, ZVP includes EVP;
but, the argument developed there is rather involved;
this is equally true for another proof of the same, proposed by
Suzuki \cite{suzuki-2001}.
A simplification of this reasoning was given in
Turinici \cite{turinici-2005},
by a technique due to
Park and Bae \cite{park-bae-1985};
note that, as a consequence of this,
ZVP $\lequi$ EVP.
It is our aim in the following to show that 
such a conclusion continues to hold  - 
under general completeness conditions - 
in the almost metric framework;
details will be given in Section 3.
The basic tool for this is a pseudometric variational principle in
Turinici \cite{turinici-2010},  
discussed in Section 2.
Finally, in Section 4 and Section 5, an application of these facts to
equilibrium points is considered.

\section{Pseudometric ordering principles}
\setcounter{equation}{0}

Let $M$ be a nonempty set.
By a {\it pseudometric} over $M$ we shall mean any map
$e:M\times M \to R_+$.
Fix such an object; which in addition is
{\it triangular} [$e(x,z)\le e(x,y)+e(y,z), \forall x,y,z\in M$] and
{\it reflexive} [$e(x,x)=0, \forall x\in M$].
Call the sequence $(x_n)$ in $M$, 
{\bf I)}
{\it strongly $e$-asymptotic} (in short: $e$-{\it strasy}) provided
[the series $\sum_n e(x_n,x_{n+1})$ converges] 
and
{\bf II)}
{\it $e$-Cauchy} when
[$\forall \de> 0, \exists n(\de)$, such that
$n(\de)\le p\le q \limpl e(x_p,x_q)\le \de$].
By the triangular property, we have
[($\forall$ sequence): $e$-strasy $\limpl$ $e$-Cauchy];
but the converse is not in general true.
Further, define an $e$-convergence structure on $X$ as: 
$x_n \etends x$ iff $e(x_n,x)\to 0$ as $n\to \oo$;
referred to as: $x$ is an $e$-limit of $(x_n)$.
The set of all these will be denoted $\lim_n(x_n)$;
when it is nonempty, we call $(x_n)$, {\it $e$-convergent}.
Note that, by the lack of symmetry, a relationship like 
[($\forall$ sequence) $d$-convergent $\limpl$ $d$-Cauchy] 
is not in general true. 
Finally, let $\vphi:M \to \Roo$ be some inf-proper function.
We consider the regularity condition
\bit
\item[(b01)]
$(e,\vphi)$ is weakly descending complete: for each $e$-strasy sequence\\ 
$(x_n)\incl \Dom(\vphi)$
with $(\vphi(x_n))$ descending there exists $x\in M$ \\
with $x_n \etends x$ and $\lim_n \vphi(x_n)\ge \vphi(x)$.
\eit
By the generic property above, it is implied by
its (stronger) counterpart
\bit
\item[(b02)]
$(e,\vphi)$ is descending complete: for each $e$-Cauchy sequence\\
$(x_n)$ in $\Dom(\vphi)$
with $(\vphi(x_n))$ descending there exists $x\in M$\\
with $x_n \etends x$ and $\lim_n \vphi(x_n)\ge \vphi(x)$.
\eit
The reciprocal inclusion also holds
[so, (b01) $\lequi$  (b02)] as it can be directly seen.

The following variational principle is our starting point.

\bprop \label{p1}
Assume that (b01)/(b02) holds. 
Then, for each $u\in \Dom(\vphi)$,  there exists 
$v=v(u)\in \Dom(\vphi)$ satisfying
(\ref{101}) (with $e$ in place of $d$) as well as

{\bf i)}
$[x\in M, e(v,x)\le \vphi(v)-\vphi(x)] \limpl \vphi(v)=\vphi(x)$
and $e(v,x)=0$.

\n
Consequently, the relations below hold

{\bf ii)}
$e(v,x)\ge \vphi(v)-\vphi(x)$, for all $x\in M$

{\bf iii)}
$e(v,x)> \vphi(v)-\vphi(x)$, for each $x\in M$ with $e(v,x)> 0$.
\eprop

The proof consists in applying 
Brezis-Browder's ordering principle  \cite{brezis-browder-1976}
to the triplet $(M_u;\le;\psi)$, where
$M_u=\{x\in M; \vphi(x)\le \vphi(u)\}$,
$(\le)$ stands for the quasi-order (a04) (with $e$ in place of $d$) and
$\psi(.)=\vphi(.)-\vphi_*$; see
Turinici \cite{turinici-2010}.

In particular, condition (b01) is retainable under
\bit
\item[(b03)]
$(e,\vphi)$ is weakly complete: for each  $e$-strasy 
sequence $(x_n)$ in $\Dom(\vphi)$ \\
there exists $x\in M$ with 
$x_n \etends x$ and $\liminf_n \vphi(x_n)\ge \vphi(x)$.
\eit
As a consequence, Proposition \ref{p1}
incorporates the variational principle in
Tataru \cite{tataru-1992};
see also
Kang and Park \cite{kang-park-1990}.

Call the pseudometric $e:M\times M\to R_+$,
an {\it almost metric} provided it is in addition
triangular and 
{\it reflexive sufficient} [$e(x,y)=0 \lequi x=y$].
In this case, Proposition \ref{p1} yields
the following practical statement.

\btheorem \label{t3}
Let the almost metric $e$
and the inf-proper function $\vphi$ be as in (b02).
Then, conclusions of Theorem \ref{t1} hold,
with $e$ in place of $d$.
\etheorem

Now, evidently, (b02) is retainable whenever
\bit
\item[(b04)]
$(e,\vphi)$ is complete: for each $e$-Cauchy
sequence $(x_n)$ in $\Dom(\vphi)$\\
there exists $x\in M$ with
$x_n \etends x$ and $\liminf_n \vphi(x_n)\ge \vphi(x)$.
\eit
If $e$ is in addition
{\it symmetric} [$e(x,y)=e(y,x), \forall x,y\in M$]
(hence, a metric over $M$), 
(b04) holds under (a02)+(a03) (modulo $e$). 
This tells us that Theorem \ref{t3} includes EVP;
it will be referred to as
the almost metric version of EVP (in short: EVPam).
The reciprocal inclusion (EVP $\limpl$ EVPam) is open;
we conjecture that the answer is positive.
Some related aspects may be found in 
Turinici \cite{turinici-1989}.

\section{Zhong variational statements}
\setcounter{equation}{0}

{\bf (A)}  Let $M$ be some nonempty set. 
Take a couple of almost metrics $d,e$ on $M$; 
we say that $e$ is $d$-{\it compatible} provided
\bit
\item[(c01)]
each $e$-Cauchy sequence is $d$-Cauchy too
\item[(c02)]
$y\mapsto e(x,y)$ is $d$-lsc, \quad for each $x\in M$.
\eit
Note that both these properties hold when $e=d$.
In fact, (c01) is trivial; 
and (c02) results from the triangular property of $d$
(see Proposition \ref{p4} for details).
Further, let $\vphi:M\to \Roo$ be an inf-proper function.
The following fact will be useful.

\blemma \label{le1}
Suppose that $e$ is $d$-compatible. Then,
[$(d,\vphi)$ = descending complete] 
implies  
[$(e,\vphi)$ = descending complete].  
\elemma

\bproof
Let $(x_n)$ be some $e$-Cauchy sequence in $\Dom(\vphi)$
with $(\vphi(x_n))$ descending.
From (c01), $(x_n)$ is $d$-Cauchy too; so, by
(b02) (modulo $d$), there exists $y\in X$ with
$x_n \dtends y$ and $\lim_n \vphi(x_n)\ge \vphi(y)$.
We claim that this is our desired point. In fact, let
$\ga> 0$ be arbitrary fixed. By the initial choice of $(x_n)$,
there exists $k=k(\ga)$ so that:
$e(x_p,x_m)\le \ga$, for each $p\ge k$ and each $m\ge p$.
Passing to limit upon $m$ one gets (via (c02))
$e(x_p,y)\le \ga$, for each $p\ge k$;
and since $\ga> 0$ was arbitrarily chosen, $x_n \etends y$.
This (and the previous relation) gives the conclusion we want.
\eproof

Now, by simply combining this with Theorem \ref{t3}, one gets
the following "relative" type variational statement 
(involving these data):

\btheorem \label{t4}
Let  the pair $(d,\vphi)$ be descending complete; 
and $e$ be $d$-compatible. 
Then, conclusions of Theorem \ref{t3} are retainable.
\etheorem

For the moment, Theorem \ref{t3} $\limpl$ Theorem \ref{t4}.
The reciprocal is also true; 
for (see above) $e=d$ is allowed here; so, 
Theorem \ref{t3} $\lequi$ Theorem \ref{t4}. 

Now, this  "relative"  variational statement 
may be viewed as an "abstract" version of ZVP. 
To explain this, we need some constructions and auxiliary facts.
\sk

{\bf (B)}
Let $b:R_+\to R_+$  be some normal function.
In particular, it is Riemann integrable on each
compact interval of $R_+$ and
\beq \label{301}
\int_p^q b(\xi)d\xi=(q-p)\int_0^1 b(p+\tau(q-p))d\tau, \ \
0\le p< q< \oo.
\eeq
Some basic facts involving the couple $(b,B)$
(where $B:R_+\to R_+$ is the one of (a06)) are being collected in

\bprop \label{p2}
The following are valid

{\bf i)} $B$ is a continuous order isomorphism of $R_+$;
hence, so is $B^{-1}$ 

{\bf ii)} $b(s)\le (B(s)-B(t))/(s-t)\le b(t), \forall t,s\in R_+, t< s$

{\bf iii)} $B$ is almost concave:
$t \mapsto [B(t+s)-B(t)]$ is decreasing on $R_+$,\ $\forall s\in R_+$

{\bf iv)} $B$ is concave:\ $B(t+\lb(s-t))\ge B(t)+\lb(B(s)-B(t))$,
for all $t,s\in R_+$ with $t< s$\ and all $\lb \in [0,1]$

{\bf v)} $B$ is sub-additive (hence $B^{-1}$ is super-additive).
\eprop

The proof is immediate, by (\ref{301})  above; so, we do not give
details. Note that the properties in {\bf iii)} and {\bf iv)}
are equivalent to each other, under {\bf i)}. This follows at once
from the (non-differential) mean value theorem in
Banta\c{s} and Turinici \cite{bantas-turinici-1994}.
\sk

{\bf (C)}
Now, let $M$ be some nonempty set; and  $d: M\times M \to R_+$, 
an almost metric over it.
Further, let $\Ga:M \to R_+$ be chosen as
\bit
\item[(c03)]
$\Ga$ is almost $d$-nonexpansive
($\Ga(x)-\Ga(y)+ d(x,y)\ge 0, \forall x,y\in M$).
\eit
Define a pseudometric $e=e(B,\Ga;d)$ over $M$ as
\bit
\item[(c04)]
$e(x,y)=B(\Ga(x)+d(x,y))-B(\Ga(x)),  \quad  x,y\in M$.
\eit
This may be viewed as an "explicit"  formula; the implicit version of
it is
\bit
\item[(c05)]
$d(x,y)=B^{-1}(B(\Ga(x))+e(x,y))-\Ga(x), \quad  x,y\in M$.
\eit
We shall establish some properties  of this map, useful in the sequel.
\sk

{\it I)}
First, the "metrical" nature of  $(x,y)\mapsto e(x,y)$ is of interest.

\bprop \label{p3}
The pseudometric $e(.,.)$ is an almost metric over $M$.
\eprop

\bproof
The reflexivity and sufficiency are clear, by Proposition \ref{p2}, i);
so, it remains to establish the triangular property.
Let $x,y,z\in M$ be arbitrary fixed. The triangular property of
$d:M\times M\to R_+$  yields (via Proposition \ref{p2}, i))
$e(x,z)\le B(\Ga(x)+d(x,y)+d(y,z))-B(\Ga(x)+d(x,y))+e(x,y)$.
On the other hand, the almost $d$-nonexpansiveness of $\Ga$ gives
$\Ga(x)+d(x,y)\ge \Ga(y)$; so (by Proposition \ref{p2}, iii))
$B(\Ga(x)+d(x,y)+d(y,z))-B(\Ga(x)+d(x,y)) \le e(y,z)$.
Combining with the previous relation yields our desired conclusion.
\eproof

{\it II)}
By definition, $e$ will be called the {\it Zhong metric} 
attached to $d$ and the couple $(B,\Ga)$.
The following  properties of  $(d,e)$ are immediate 
(via Proposition \ref{p2}):

\blemma \label{le2}
Under the prescribed conventions,

{\bf vi)}
$b(\Ga(x)+d(x,y))d(x,y)\le e(x,y)\le b(\Ga(x))d(x,y)$,\ for all $x,y\in M$

{\bf vii)}
$e(x,y)\le B(d(x,y)), \forall x,y\in M$;
hence  $x_n\dtends x$ implies $x_n\etends x$.
\elemma

{\it III)}  A basic property of $e(.,.)$ to be checked is
$d$-compatibility.

\bprop \label{p4}
The Zhong metric $e$ is $d$-compatible (cf. (c01)+(c02)).
\eprop

\bproof
We firstly check (c02); which may be written as
$$
\mbox{
[$e(x,y_n)\le \lb, \forall n$]\ and $y_n\dtends y$ \ \ imply  $e(x,y)\le \lb$.
}
$$
So, let $x$, $(y_n)$, $\lb$ and $y$ be as in the premise of this relation.
By Lemma \ref{le2}, we have $y_n\etends y$ as $n\to \oo$.
Moreover (as $e$ is triangular)
$e(x,y)\le e(x,y_n)+e(y_n,y)\le \lb+e(x_n,y)$, for all $n$.
It will suffice passing to limit as $n\to \oo$
to get the desired conclusion.
Further, we claim that (c01) holds too, in the sense:
[(for each sequence)\ \ $d$-Cauchy $\lequi$  $e$-Cauchy].
The left to right implication is clear, via Lemma \ref{le2}.
For the right to left one, assume that $(x_n)$ is an $e$-Cauchy
sequence in $M$. In particular (by the triangular property)
$e(x_i,x_j)\le \mu$, for all $(i,j)$ with $i\le j$
and some $\mu\ge 0$.
This, along with (c05), yields
$d(x_0,x_i)=B^{-1}(B(\Ga(x_0))+e(x_0,x_i))-\Ga(x_0)\le$
$B^{-1}(B(\Ga(x_0))+\mu)-\Ga(x_0)$, $\forall i\ge 0$;
wherefrom (cf. (c03))
$\Ga(x_i)\le \Ga(x_0)+d(x_0,x_i)\le B^{-1}(B(\Ga(x_0))+\mu)$
(hence $B(\Ga(x_i))\le B(\Ga(x_0)+\mu$), for all $i\ge 0$.
Putting these facts together yields (again via (c05))
$\Ga(x_i)+d(x_i,x_j)=B^{-1}(B(\Ga(x_i))+e(x_i,x_j))\le$ 
$\nu:=B^{-1}(B(\Ga(x_0))+2\mu)$, for all $(i,j)$ with $i\le j$.
And this, via Lemma \ref{le2}, gives (for the same pairs $(i,j)$)
$e(x_i,x_j)\ge b(\Ga(x_i)+d(x_i,x_j))d(x_i,x_j)\ge
b(\nu)d(x_i,x_j)$.
But then, the $d$-Cauchy property of $(x_n)$ is clear; 
and the proof is complete.
\eproof

{\bf (D)}
We are now in position to make precise our initial claim.
Let the  almost metric $d$ and the inf-proper function $\vphi$ be as in
(b02) (modulo $d$).  Further, take a normal function $b:R_+\to R_+$;
as well as an almost $d$-nonexpansive map $\Ga:M\to R_+$.
Finally, put $e=e(B,\Ga;d)$ (the Zhong metric introduced by (c04)/(c05)).

\btheorem \label{t5}
Let the conditions above be admitted. Then

{\bf viii)}
For each  $u\in \Dom(\vphi)$
there exists $v=v(u)\in \Dom(\vphi)$ with
\beq \label{302}
b(\Ga(u)+d(u,v))d(u,v)\le e(u,v)\le \vphi(u)-\vphi(v)
\eeq
\beq \label{303}
b(\Ga(v))d(v,x)\ge e(v,x)> \vphi(v)-\vphi(x), \quad
\forall x\in M\sm \{v\}
\eeq

{\bf ix)}
For each $u\in \Dom(\vphi)$,  $\rho> 0$  with
$\vphi(u)-\vphi_*\le B(\Ga(u)+\rho)-B(\Ga(u))$
the above evaluation (\ref{302}) gives
\beq \label{304}
\mbox{
$d(u,v)\le \rho$; \quad hence $\Ga(v)\le \Ga(u)+\rho$
}
\eeq
\beq \label{305}
\mbox{
$b(\Ga(u)+\rho)d(u,v)\le \vphi(u)-\vphi(v)$ \ \
(hence $\vphi(u)\ge \vphi(v)$).
}
\eeq
\etheorem

\bproof
By Proposition \ref{p3}, $e$ is an almost metric over $M$;
and, by Proposition \ref{p4}, it is $d$-compatible. 
Hence, Theorem \ref{t4} applies to such data.
In this case, (\ref{302})+(\ref{303}) are clear via Lemma \ref{le2}.
Moreover, if $u\in \Dom(\vphi)$ is taken as in
the premise of {\bf ix)}, then (cf. (\ref{302}))
$e(u,v)\le \vphi(u)-\vphi(v)\le \vphi(u)-\vphi_*$; wherefrom (by (c05))
$d(u,v)\le B^{-1}(B(\Ga(u))+\vphi(u)-\vphi_*)-\Ga(u)\le \rho$;
and (\ref{304})+(\ref{305}) follow as well.
\eproof

So far, 
Theorem \ref{t3} $\limpl$ Theorem \ref{t4} $\limpl$ Theorem \ref{t5}.
In addition, Theorem \ref{t5} $\limpl$ Theorem \ref{t3};
just take $b=1$ (hence $B$=identity, $e=d$).
Summing up, these three variational principles are mutually equivalent.
On the other hand, Theorem \ref{t5} may be also viewed as an
extended (modulo $\Ga$) version of ZVP.
For, if $d$ is symmetric (hence a (standard) metric),  (c03) becomes
\bit
\item[(c06)]
$|\Ga(x)-\Ga(y)|\le d(x,y)$,\ for all $x,y\in M$\ \ 
($\Ga$ is $d$-nonexpansive).
\eit
In addition, the choice
\bit
\item[(c07)]
$\Ga(x)=d(a,x), \ x\in M$, \quad for some $a\in M$
\eit
is in agreement with it; hence the claim.
For this reason, Theorem \ref{t5} will be referred to as
the almost metric version of ZVP (in short: ZVPam).
This inclusion is technically strict; because
the  conclusions involving the middle terms in
(\ref{302})+(\ref{303}) cannot be obtained in the way described by
Zhong \cite{zhong-1997}.
Some related aspects were delineated in
Ray and Walker \cite{ray-walker-1982};
see also
Suzuki \cite{suzuki-2001}.

\section{Application (equilibrium points)}
\setcounter{equation}{0}

Let $M$ be some nonempty set;
and $e:M\times M\to R_+$ be an almost metric over it.
Any (extended) function
$G:M\times M\to R\cup \{-\oo,\oo\}$
will be referred to as a 
{\it relative generalized pseudometric} on $M$.
Given such an object, 
we say that $v\in M$ is an {\it equilibrium} point of 
it, when $G(v,x)\ge 0$, $\forall x\in M$.
Note that, if $G=F+e$, where
$F:M\times M\to R\cup \{-\oo,\oo\}$
is another relative generalized pseudometric on $M$,
this definition becomes
$e(v,x)\ge -F(v,x)$, $\forall x\in M$;
which, under the choice (for some $\vphi:M\to \Roo$)
\bit
\item[(d01)]
$F(x,y)=\vphi(y)-\vphi(x), \ \ x,y\in M$ \quad (where $\oo-\oo=0$)
\eit
tells us that  the variational property  of $v$ is "close" 
to the one in Theorem \ref{t3}.
So, existence of such points is deductible from the quoted result;
to do this, one may proceed as follows.
Assume that the relative generalized pseudometric $F$ is
{\it triangular}
[$F(x,z)\le F(x,y)+F(y,z)$, whenever the right member exists] and
{\it reflexive}
[$F(x,x)=0$, $\forall x\in M$].
Define  the (extended) function
$\mu(x)=\sup\{-F(x,y); y\in M\}, x\in M$;
clearly, $\mu(M)\incl \R+oo$. 
The alternative $\mu(M)=\{\oo\}$ cannot
be excluded; to avoid this, assume
\bit
\item[(d02)]
$\mu$ is proper ($\Dom(\mu):=\{x\in M; \mu(x)< \oo\} \ne \es$).
\eit
For the arbitrary fixed $u\in \Dom(\mu)$ put $F_u(.)=F(u,.)$.
We have by definition
\beq \label{401}
F_u(u)=0;\  F_u^*:=\inf\{F_u(x); x\in M\}=-\mu(u)>-\oo;
\eeq
so,  $F_u$ is inf-proper [and we say: $F$ is semi inf-proper].
Further, let $d$ be another almost metric on $M$ with
\bit
\item[(d03)]
$(d,F)$ is semi descending complete:\\
$(d,F_u)$ is descending complete,  for each $u\in \Dom(\mu)$.
\eit

\btheorem \label{t6}
Let (d02)+(d03) hold; and $e$ be $d$-compatible.
Then, for each $u\in \Dom(\mu)$ there exists $v=v(u)$ in $M$ such that

{\bf i)}
$e(u,v)\le -F(u,v)\le \mu(u)(< \oo)$

{\bf ii)}
$e(v,x)> -F(v,x)$, for all $x\in M\sm \{v\}$.

\n
Hence, in particular, $v$ is an equilibrium point for $G:=F+e$.
\etheorem

\bproof
From Theorem \ref{t4} it follows that, for the
starting $u\in \Dom(\mu)$ (hence $u\in \Dom(F_u)$) there
must be another point $v\in \Dom(F_u)$ with the properties
{\bf I)} $e(u,v)\le F_u(u)-F_u(v)$ 
and
{\bf II)} $e(v,x)> F_u(v)-F_u(x), \forall x\in M\sm \{v\}$.
The former of these is just {\bf i)},
by the reflexivity of $F$.
And the latter yields {\bf ii)}; 
for (by the triangular property)
$F(u,v)-F(u,x)\ge F(u,v)-(F(u,v)+F(v,x))=-F(v,x)$.
\eproof

Now, a basic particular choice of $e(.,.)$ is
related to the constructions in Section 3. Precisely, let
the function $b:R_+\to R_+$ be normal;
and $\Ga:M\to R_+$ be almost $d$-nonexpansive.
Let $e=e(B,\Ga;d)$ stand for the 
Zhong metric given by (c04)/(c05).
By Theorem \ref{t5}, we then have

\btheorem \label{t7}
Let (d02)+(d03) hold. 
Then, for each $u\in \Dom(\mu)$ there exists $v=v(u)$ in $M$ such that

{\bf iii)}
$b(\Ga(u)+d(u,v))d(u,v)\le e(u,v)\le -F(u,v)\le \mu(u)$

{\bf iv)}
$b(\Ga(v))d(v,x)\ge e(v,x)> -F(v,x), \forall x\in M\sm \{v\}$.

\n
Hence, in particular, $v$ is an equilibrium point for 
$G(x,y)=F(x,y)+b(\Ga(x))d(x,y)$, $x,y\in M$.
Moreover, $u\in \Dom(\mu)$ whenever
\bit
\item[(d04)]
$\mu(u)\le B(\Ga(u)+\rho)-B(\Ga(u))$, for some $\rho> 0$;
\eit
and then (as  $F_u(u)-F_u^*=\mu(u)$),
{\bf iii)} gives (\ref{304}) and

{\bf v)}
$b(\Ga(u)+\rho)d(u,v)\le -F(u,v)$ (hence $F(u,v)\le 0$).
\etheorem

Some remarks are in order. Let $\vphi:M\to \Roo$ be some inf-proper
function. The relative (generalized) pseudometric 
$F$ over $M$ given as in (d01)
is reflexive, triangular and fulfills (d02); because
$\mu(.)=\vphi(.)-\vphi_*$ (hence $\Dom(\mu)=\Dom(\vphi)$).
In addition, as
$F_u(.)=\vphi(.)-\vphi(u)$, $u\in \Dom(\vphi)$,
(d03) is identical with (b02) (modulo $d$). 
Putting these together, it follows that
Theorems \ref{t6} and \ref{t7}  include
Theorems \ref{t4} and \ref{t5} respectively.
The reciprocal inclusions are also true, by the very argument above;
so that
Theorem \ref{t6} $\lequi$ Theorem \ref{t4}  and
Theorem \ref{t7} $\lequi$ Theorem \ref{t5}.
In particular, when $\Ga$ is taken as in (c07),
Theorem \ref{t7} yields the main result in
Zhu, Zhong and Cho \cite{zhu-zhong-cho-2000};
see also
Bao and Khanh \cite{bao-khanh-2003}.

\section{The BKP approach}
\setcounter{equation}{0}

Let $(M,d)$ be a complete metric space.
By a {\it relative pseudometric} over $M$ we mean 
any map $g:M\times M\to R$. Given such an object, 
remember that $v\in M$ is an equilibrium point of it 
when $g(v,x)\ge 0$, $\forall x\in M$.
Note that, if $g=f+d$, where $f:M\times M\to R$ 
is another relative pseudometric on $M$, this writes 
$d(v,x)\ge -f(v,x)$, $\forall x\in M$; 
so that, under the choice (d01) of $f$ (where $\vphi:M\to R$),
the variational property of  $v$ is "close" to the one in EVP.
The folllowing 2005 result in the area due to
Bianchi, Kassay and Pini \cite{bianchi-kassay-pini-2005}
(in short: BKP) is available.

\btheorem \label{t8}
Suppose that $f$ is triangular reflexive and
\bit
\item[(e01)]
$f(a,.)$ is bounded from below and lsc, for each $a\in M$.
\eit
Then, for each $u\in M$, there exists $v=v(u)\in M$ such that

{\bf i)}
$d(u,v)\le -f(u,v)$;\ \ 
{\bf ii)}
$d(v,x)>  -f(v,x)$,\ for all\ $x\in M\sm \{v\}$.

\n
Hence, in particular, $v$ is an equilibrium point for $g:=f+d$.
\etheorem

Note that this result is obtainable from Theorem \ref{t6}
by simply taking $d=e$.
On the other hand, under the same choice (d01) for $f$,
(e01) becomes
\bit
\item[(e02)]
$\vphi$ is bounded from below and lsc;
\eit
and Theorem \ref{t8} is just EVP.
So, we may ask whether this extension is effective.
The answer is negative;
i.e., Theorem \ref{t8} is deductible from 
(hence equivalent with) EVP. This will follow from

\bproof  ({\bf Theorem \ref{t8}})
Define a new function $h:M\to R$ as 
$h(x)=f(u,x)$, $x\in M$.
From (e01), EVP is applicable to $(M,d)$ and $h$;
wherefrom, for the starting $u\in M$ there exists 
$v\in M$ with
{\bf I)} $d(u,v)\le h(u)-h(v)$,
{\bf II)} $d(v,x)> h(v)-h(x),\ \ \forall x\in M\sm \{v\}$.
The former of these gives  {\bf i)}, in view of
$h(u)=0$.
And the latter one gives {\bf ii)};
because (from the triangular property)
$h(v)-h(x)\ge -f(v,x)$, for all such $x$.
Hence the conclusion.
\eproof

This argument 
(taken from the 2003 paper due to
Bao and Khanh \cite{bao-khanh-2003})
tells us that Theorem \ref{t8} 
is just a formal extension of EVP.
This is also true for the 1993 statement in the area due to
Oettli and Thera \cite{oettli-thera-1993}.
In fact, the whole reasoning developed in 
\cite{bianchi-kassay-pini-2005}
for proving Theorem \ref{t8}
is, practically,  identical with the one of this last paper.
Further aspects may be found in 
Lin and Du \cite{lin-du-2006}.


\end{document}